\documentclass{elsart}

\usepackage{natbib}
\usepackage{graphicx}
\usepackage{amssymb, amsmath}

\begin{document}

\begin{frontmatter}

\title{Computing all roots of the likelihood equations
of seemingly unrelated regressions}

\author{Mathias Drton\thanksref{support}}

\address{Department of Statistics, The University of Chicago, Chicago,
  IL  60637, U.S.A.}
\thanks[support]{Supported by the University of Washington Royalty Research
Fund Grant No.\ 65-3010.} 

\begin{abstract}
Seemingly unrelated regressions are statistical regression
models based on the Gaussian distribution.  They are popular in
econometrics but also arise in graphical modeling of 
multivariate dependencies.  In maximum likelihood estimation, the
parameters of the model are estimated by maximizing the likelihood
function, which maps the parameters to the likelihood of observing the
given data.  By transforming this optimization problem into a
polynomial optimization problem, it was recently shown that the
likelihood function of a simple bivariate seemingly unrelated
regressions model may have several stationary points. Thus local
maxima may complicate maximum likelihood estimation. In this paper, we
study several more complicated seemingly unrelated regression models,
and show how all stationary points of the likelihood function can be
computed using algebraic geometry.
\end{abstract}


\end{frontmatter}

%% MACROS %%%%%%%%%%%%%%%%%%%%%%%%%%%%%%%%%
\newcommand{\RRR}{\mathbb{R}}
\newcommand{\QQQ}{\mathbb{Q}}
\newcommand{\CCC}{\mathbb{C}}
\newcommand{\ND}{\mathcal{N}}
\newcommand{\CRfam}{\mathcal{C_R}}
\newcommand{\NSUR}{\mathbf{N}(\CRfam)}
\newcommand{\BCR}{\mathbb{B}(\CRfam)}
\newcommand{\ellprof}{ \ell_{\mathrm{prof}}}
%%%%%%%%%%%%%%%%%%%%%%%%%%%%%%%%%%%%%%%%%%%

\section{Introduction}
\label{sec:intro}

Seemingly unrelated regressions (SUR) are multivariate regression
models with correlated response (or dependent) variables that follow a
joint Gaussian distribution. Usually different regressions contain
different covariates (or independent variables) and seem
``unrelated.''  However, due to the correlated response variables the
regressions are only ``seemingly unrelated'' and contain valuable
information about each other \citep{zellner:62}. SUR play ``a central
role in contemporary econometrics'' \cite[p.~323]{goldberger:91} but
also appear in other contexts
\citep{rochon:jrss,rochon:bio,verbyla:1988}.  Moreover, SUR arise in
the context of Gaussian graphical models
(\citealp[\S5]{amp:cgmp}; \citealp[\S8.5]{rich:annals}).  

The parameters of a SUR model can be estimated efficiently, i.e.\ with
small variance, by maximizing the likelihood function, which maps the
parameters to the likelihood of observing the given data.
\cite{oberhofer:74} and \cite{telser:64} give two popular algorithms
for this maximization.  In general, however, these algorithms will not
globally maximize the likelihood function, which indeed may be
multimodal; a fact neglected in the literature \cite[\S6]{drton:2004}.
\cite{drton:2004} demonstrated the possibility of multimodality in a
study of a bivariate SUR model that may have a likelihood function
with five stationary points. In this paper, we use algebraic geometry
to apply the approach of \cite{drton:2004} to more general SUR models.
In Sections \ref{sec:sur} and \ref{sec:mlepoly} we give an
introduction to SUR and show how maximum likelihood estimation can be
performed by solving a polynomial optimization problem, opening the
door for tools from algebraic geometry.  With
these tools, we first revisit the work by \cite{drton:2004}, see
Section \ref{sec:revisit}, and then obtain new results on more general
SUR models (Section \ref{sec:groebner}).  In particular, we identify
examples of SUR models, for which all stationary points of the
likelihood function can be computed.

\section{Seemingly unrelated regressions}
\label{sec:sur}

In SUR a family of response variables, indexed by a finite set $R$, is
stochastically modeled using a family of covariates, indexed by a finite
set $C$.  All response variables and all covariates are observed on a
finite set of subjects $N$. We denote the cardinalities of the three sets
also by $R$, $C$ and $N$, respectively. The
observations can be represented by two matrices $X$ and $Y$.  The matrix
$Y=(Y_{rm})\in\RRR^{R\times N}$ has the $(r,m)$-entry equal to the
observation of response variable $r\in R$ on subject $m\in N$, and the
matrix $X=(X_{cm})\in \RRR^{C\times N}$ has the $(c,m)$-entry equal to the
observation of covariate $c\in C$ on subject $m\in N$.  For $c\in C$ and
$r\in R$, $X_c\in\RRR^N$ and $Y_r\in \RRR^N$ denote the $c$-th and $r$-th
row of $X$ and $Y$, respectively.  Similarly, $X^m$ and $Y^m$, $m\in N$,
denote the $m$-th column of $X$ and $Y$, respectively. Clearly, $X_c$ and
$Y_r$ comprise all observations of the $c$-th covariate and the $r$-th
response variable; $X^m$ and $Y^m$ comprise all covariate and response
variable observations on the $m$-th subject.

In this regression setting, the matrix $X$ is assumed to be
deterministic and fixed but the matrix $Y$ is modeled to follow a
multivariate normal distribution, where the mean vector of $Y_r$, $r\in R$,
is a linear combination of some $X_c$, $c\in C_r\subseteq C$,
\begin{equation}
  \label{eq:betadef}
  \mathrm{E}[Y_r] = \sum_{c\in C_r} \beta_{rc} X_c\in \RRR^N, \quad r\in R. 
\end{equation}
Here $(C_r\mid r\in R)$ is a fixed family of subsets of
$C$ indexing the covariates involved in each one of the $R$ regressions. The
weights $\beta_{rc}$ in (\ref{eq:betadef}) are called {\em regression
  coefficients\/}.  Setting $\beta_{rc}=0$ if $c\not\in C_r$, 
we can define a matrix of regression coefficients
$B=(\beta_{rc})\in \RRR^{R\times C}$. The random vectors $Y^m$, $m\in
N$, are assumed to be independent with common positive definite covariance
matrix  
\begin{equation}
  \label{eq:sigmadef}
  \mathrm{Var}[Y^m]=\Sigma\in \RRR^{R\times   R}, \quad m\in N.
\end{equation}
Letting
\begin{equation}
  \label{eq:CRdef}
  \CRfam=\cup (\{r\}\times C_r \mid r\in R) \subseteq R\times C,
\end{equation}
the {\em seemingly unrelated regressions model\/} is the family
of normal distributions 
\begin{equation}
  \label{eq:surdef}
  \NSUR = \big( \ND_{R\times N}(BX,\Sigma\otimes I_N) \,\big|\,
  (B,\Sigma)\in \BCR\times \mathbb{P} \big). 
\end{equation}
Here $\ND_{R\times N}$ is the multivariate normal distribution on
$\RRR^{R\times N}$; $I_N$ is the $N\times N$ identity matrix; $\otimes$ is
the Kronecker product; $B$ and $\Sigma$ are the mean and the variance
parameters; and {\em the parameter space\/} $\BCR\times \mathbb{P}$ is the
Cartesian product of the linear space
\begin{equation}
  \label{eq:bspacedef}
  \BCR = \big\{ B\in \RRR^{R\times C} \,\big|\,
  B=(\beta_{rc}),\; \beta_{rc}=0 \;\forall (r,c)\not\in \CRfam \big\}
\end{equation}
and the cone $\mathbb{P}$ of all positive definite real $R\times R$ matrices.
The response matrix $Y$ is then an observation from some (unknown) distribution
in the model,
\begin{equation*}
  Y \sim \ND(BX, \Sigma \otimes I_N), \quad (B,\Sigma)\in
  \BCR\times \mathbb{P}. 
\end{equation*}
% {\bf Write out regression equations}
If $N\ge R+C$ and $X$ is a matrix of full rank, then with probability
one the $(R+C)\times N$ matrix obtained by stacking $X$ and $Y$ has
full rank,  
\begin{equation}
  \label{eq:fullrank}
  \text{rank} 
  \begin{pmatrix}
    Y\\X
  \end{pmatrix} = R+C.
\end{equation}
We assume (\ref{eq:fullrank}) to hold throughout the paper.

\section{Maximum likelihood estimation by polynomial optimization}
\label{sec:mlepoly}

The {\em probability density function\/} $f_{(B,\Sigma)}: \RRR^{R\times N} 
\to (0,\infty)$ of the distribution
$\mathcal{N}(BX,\Sigma\otimes I_n)$ can be written as 
\begin{equation*}
  f_{(B,\Sigma)}(Y) = \frac{1}{\sqrt{(2\pi)^{RN}
      |\Sigma|^N}}\exp\left\{-\frac{1}{2}
    \mathrm{tr}\big[\Sigma^{-1}(Y-BX)(Y-BX)'\big] \right\}. 
\end{equation*}
For data $Y$, the {\em likelihood function\/} $L: \BCR\times \mathbb{P}
\to (0,\infty)$ of the model $\NSUR$ is 
defined as
\begin{equation*}
  L(B,\Sigma) = f_{(B,\Sigma)}(Y). 
\end{equation*}
In maximum likelihood estimation the parameters $(B,\Sigma)$ are
estimated by
\begin{equation}
  \label{eq:maxL}
  (\hat B, \hat\Sigma) = \arg\max \{ L(B,\Sigma) \mid (B,\Sigma)\in
  \BCR\times \mathbb{P} \} .
\end{equation}
It follows from (\ref{eq:fullrank}) that the maximum of the
likelihood function exists. 

We can parameterize $\BCR$ by mapping a vector
\begin{equation*}
  \beta=(\beta_{rc}\mid (r,c)\in \CRfam) \in \RRR^{\CRfam},
\end{equation*}
to the matrix $B(\beta)\in \BCR$ with entry $B(\beta)_{rc}=\beta_{rc}$ if
$(r,c)\in\CRfam$ and $B(\beta)_{rc}=0$ otherwise.  
Define $\ell: \RRR^{\CRfam}\times \mathbb{P} \to \RRR$ by 
  \begin{equation}
    \label{eq:elldef}
    \begin{split}
      \ell(\beta,\Sigma) &= \log L(B(\beta),\Sigma) 
      \\ &
      \propto
      -\frac{N}{2}\log |\Sigma| - \frac{1}{2}
      \mathrm{tr}\big[\Sigma^{-1}
      \big(Y-B(\beta)X\big)\big(Y-B(\beta)X\big)'\big].    
    \end{split}
  \end{equation}
Clearly we can solve (\ref{eq:maxL}) by finding
\begin{equation}
  \label{eq:maxell}
  (\hat \beta, \hat\Sigma) = \arg\max \{ \ell(\beta,\Sigma) \mid
  (\beta,\Sigma)\in \RRR^{\CRfam}\times \mathbb{P} \},
\end{equation}
and setting $\hat B=B(\hat \beta)$.  
The standard approach to solve (\ref{eq:maxell}) is to solve the 
{\em likelihood equations\/}  
\begin{equation}
  \label{eq:likeqn}
  \left(\frac{\partial \ell(\beta,\Sigma)}{\partial \beta},
  \frac{\partial \ell(\beta,\Sigma)}{\partial \Sigma}\right) =0. 
\end{equation}
It can be shown that (\ref{eq:likeqn}) holds if and only if
\begin{equation}
  \label{eq:sigmamax}
  \Sigma = \frac{1}{N} \big(Y-B(\beta)X\big)\big(Y-B(\beta)X\big)' 
\end{equation}
and 
\begin{equation}
  \label{eq:betamax}
  \beta = \left[ A'(XX'\otimes \Sigma^{-1}) A \right]^{-1}
  A' \mathrm{vec}(\Sigma^{-1} YX'),
\end{equation}
where $A$ is a matrix of zeroes and ones that satisfies
$\mathrm{vec}(B(\beta))=A\beta$.  In fact, each column of $A$ has precisely
one entry equal to one and the remaining entries equal to zero.
\cite{oberhofer:74} show how one solution to the likelihood equations
can be obtained by alternating between solving
(\ref{eq:sigmamax}) for fixed $\beta$ and solving (\ref{eq:betamax})
for fixed $\Sigma$.
Here, we take a different approach that, for certain SUR models,
allows us to compute all solutions to the likelihood equations.

From (\ref{eq:fullrank}) and (\ref{eq:elldef}), it follows that for
fixed $\beta\in\RRR^{\CRfam}$ the function $\ell_\beta: \Sigma \mapsto
\ell(\beta,\Sigma)$ is strictly concave with maximizer
(\ref{eq:sigmamax}).   Thus the {\em profile 
log-likelihood function\/} $\ell_{\mathrm{prof}}: \RRR^{\CRfam} \to \RRR$
defined as 
\begin{equation}
  \label{eq:profell}
  \ellprof (\beta) = \max\{\ell(\beta,\Sigma) \mid
  \Sigma\in\mathbb{P}\}
\end{equation}
takes on the form 
\begin{equation}
  \label{eq:profellsimple}
  \ellprof (\beta) \propto -\frac{N}{2}\log \big| \frac{1}{N}
  \big(Y-B(\beta) X\big)\big(Y-B(\beta)X\big)' \big| - \frac{RN}{2}.
\end{equation}
By the strict con-cavity of $\ell_\beta$, $(\beta,\Sigma)$ is a stationary
point of $\ell(\beta, \Sigma)$ if and only if $\beta$ is a
stationary point of $\ellprof(\beta)$ and $\Sigma$ satisfies
(\ref{eq:sigmamax}); compare \citet[Lemma 1]{drton:2004}.  
The same holds for
\begin{equation}
  \label{eq:defG}
  G(\beta)= \left| \big(Y-B(\beta)X\big)\big(Y-B(\beta)X\big)' \right|,
\end{equation}
which conveniently is a polynomial in $\beta$.  Thus we can solve
(\ref{eq:maxell}) by using (\ref{eq:sigmamax}) and solving the
unconstrained polynomial program 
\begin{equation}
  \label{eq:minG}
  \hat\beta = \arg\min \{ G(\beta) \mid \beta\in\RRR^{\CRfam} \}.
\end{equation}

%% \section{Maximum likelihood ideals}
%% \label{sec:alg}

We try to solve (\ref{eq:minG}) by
computing the stationary points of $G$, i.e. by solving the equations
\begin{equation}
  \label{eq:defgrc}
  g_{rc} = \frac{\partial G(\beta)}{\partial\beta_{rc}} = 0, \quad
  (r,c)\in\CRfam.
\end{equation} 
In practice the observations
$Y$ and $X$ are available only in finite accuracy and the partial
derivatives $g_{rc}$, $(r,c)\in\CRfam$, are elements of the ring $\QQQ[\beta]$
of polynomials in $\beta$ with rational
coefficients. 
In an algebraic approach to solving polynomial equations 
\citep{coxlittleoshea:1997,coxlittleoshea:1998,sturmfels:2002} we
allow the indeterminants in the polynomial equation system
(\ref{eq:defgrc}) to be complex, i.e.\ 
$\beta\in\CCC^{\CRfam}$, where $\CCC$ is the field of complex numbers.
We define the {\em maximum likelihood ideal\/} $I_G$ to be the ideal
of $\QQQ[\beta]$ that is generated by the partial derivatives
$g_{rc}$, $(r,c)\in\CRfam$, i.e.\
\begin{equation}
  \label{eq:defIG}
  I_{G} = \langle g_{rc} \mid (r,c)\in\CRfam \rangle;
\end{equation}
compare \citet[\S8.4]{sturmfels:2002} who defines maximum likelihood ideals
in a different statistical context.  Software like {\tt Macaulay
  2}\footnote{\tt http://www.math.uiuc.edu/Macaulay2/} and {\tt Singular}
\citep{singular} permits us to check whether $I_G$ is a zero-dimensional
ideal.  If $\dim(I_G)=0$, then the variety $V_\CCC(I_G)$, i.e.\ the set of
common complex zeroes of the partial derivatives $g_{rc}$, is a finite set
and all its elements can be computed using, for example, 
{\tt Singular} or also {\tt PHCpack}\footnote{\tt
  http://www.math.uic.edu/\~{}jan/}. The real points   
$V_\RRR(I_G)=V_\CCC(I_G)\cap \RRR^\CRfam$ can then be identified and yield
the stationary points of $G$.

\section{Revisiting the multimodal bivariate seemingly unrelated
  regressions with two covariates}
\label{sec:revisit}

\cite{drton:2004} study a SUR model with two response variables and two
covariates, in which response variable 1 is regressed only on covariate 1,
and response variable 2 only on covariate 2.  Hence, $R=\{1,2\}$, 
$C=\{1,2\}$, $C_1=\{1\}$, and $C_2= \{2\}$.  Therefore, $\CRfam=\{(1,1), 
(2,2)\}$, and $B\in\BCR$ if $B$ is of the form 
\begin{equation*}
  B= \begin{pmatrix}
    \beta_{11} & 0\\
    0 & \beta_{22}
  \end{pmatrix} \in \RRR^{2\times 2}.
\end{equation*}
Using {\tt Singular} and the data
in \citet[Table 1]{drton:2004},  we can solve (\ref{eq:minG}) as shown in
Table \ref{tab:singularcode}.
\begin{table}[tb]
\begin{center}
\hrule
\vspace{0.05cm}
\hrule
{ 
\medskip
\begin{verbatim}
> ring R=0,(b(1..2)), lp;
> matrix X[2][8] = 188,22,-46,77,-103,74,83,101,   
.                  55,-216,116,-30,131,195,-311,-239;
> matrix Y[2][8] = 234,-5,6,182,-193,278,62,-68,   
.                  497,-326,266,-3,93,558,-584,-224;
> matrix B[2][2] = b(1),0, 0,b(2);
> poly G = det((Y-B*X)*transpose(Y-B*X));
> ideal IG =jacob(G);  
> ideal J = groebner(IG);
> dim(J); vdim(J);
0
5
> LIB "solve.lib"; solve(J,6);
[1]:
   [1]:  0.778796
   [2]:  1.538029
[2]:
   [1]:  1.622609
   [2]:  2.034745
[3]:
   [1]:  (1.480687-i*1.547274)
   [2]:  (2.16845+i*0.765283)
[4]:
   [1]:  (1.480687+i*1.547274)
   [2]:  (2.16845-i*0.765283)
[5]:
   [1]:  2.764418
   [2]:  2.504006
\end{verbatim} 
}
\medskip
\hrule
\end{center}
\medskip
\caption{{\tt Singular}-session for the model in
\cite{drton:2004}.\label{tab:singularcode} }
\end{table}

As computed by {\tt dim} and {\tt vdim}, the maximum likelihood ideal
$I_G={\tt IG}$ is zero-dimensional and of degree five. The five points
in the variety $V_\CCC(I_G)$ are computed 
by {\tt solve}, which lists $\beta_{11}={\tt b(1)}$ as first component
and $\beta_{22}={\tt b(2)}$ as second component. There are three real
points in $V_\RRR(I_G)$, which yield the 
stationary points of the likelihood function of the model $\NSUR$. 
Note that we confirm the values
stated in \citet[Table 2 with $\beta_{11}=\beta_1$ and 
$\beta_{22}=\beta_2$]{drton:2004}. 
The Gr\"obner basis computed by the command {\tt 
groebner(IG)} has two elements that are (i)
a quintic in $\beta_{22}={\tt b(2)}$ and (ii) a sum of a linear
function in $\beta_{11}={\tt b(1)}$ and a quartic in $\beta_{22}={\tt
b(2)}$. Thus it follows immediately that the stationary
points of $G$ can be found from solving a quintic \cite[cf.][Thm.\
2]{drton:2004}.

\section{Dimensions and degrees of maximum likelihood ideals}
\label{sec:groebner}

\subsection{Seemingly unrelated regressions} 
\label{sec:dimdegsur}

The algebraic approach can also be applied to more
general models.  Here we focus on SUR models $\NSUR$ for which
$(C_r\mid r\in R)$ consists of disjoint sets; in other models 
inclusion relations among the sets $C_r$ may be
exploited \citep[cf.][]{ap:ims}.  
More precisely, we consider models $\NSUR$ in
which $r_1< r_2$,  $r_1,r_2\in R$, implies that $c_1< c_2$ for all 
$c_1\in C_{r_1}$ and $c_2\in C_{r_2}$.  Then $\BCR$
is a linear space of block-diagonal matrices.

\begin{figure}[p] 
\centering
\begin{minipage}[t]{6.5cm}
\begin{center} 
\includegraphics[width=6.5cm]{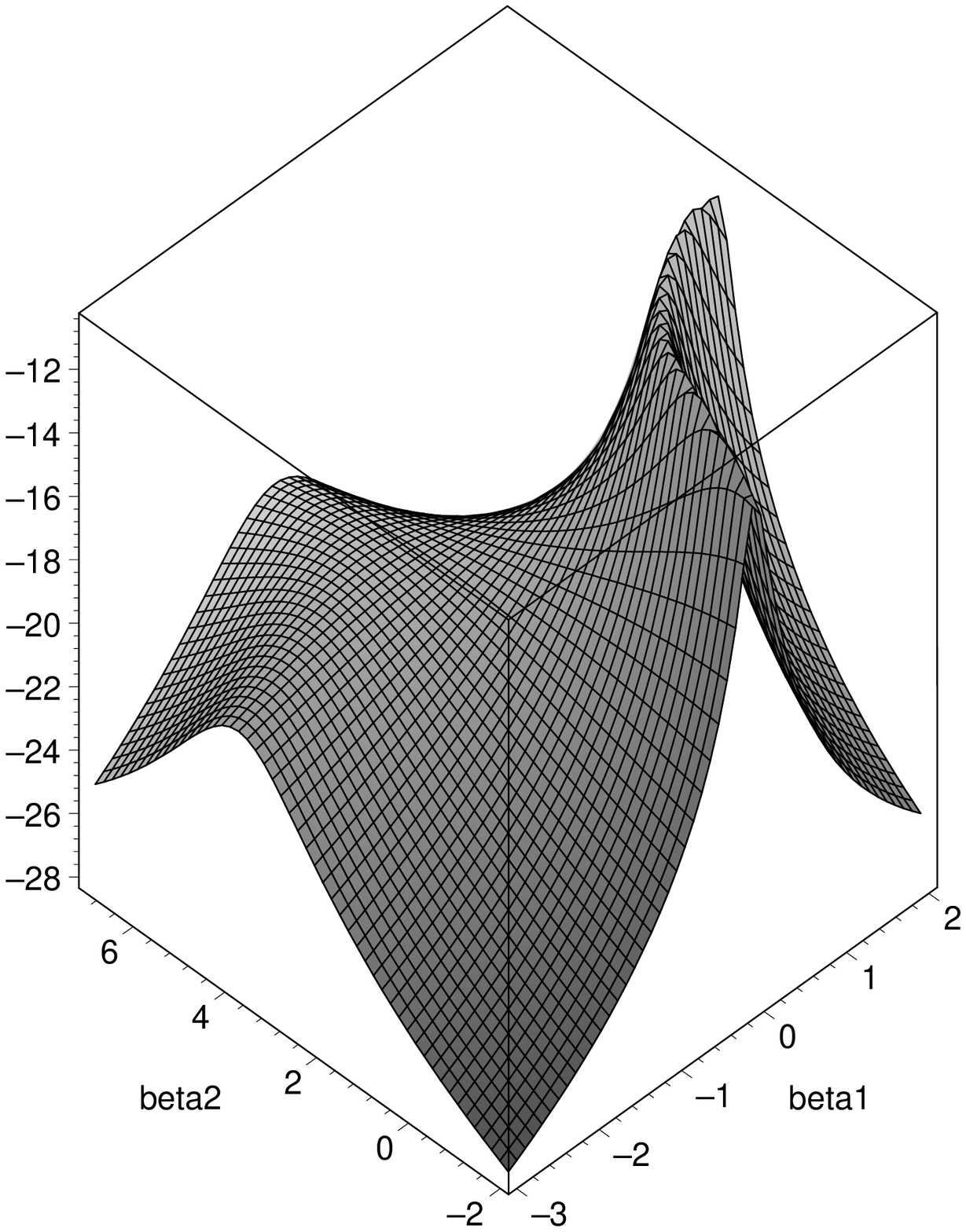} 
\end{center} 
\end{minipage}
\hfill 
\begin{minipage}[t]{6.5cm}
\begin{center} 
\includegraphics[width=6cm]{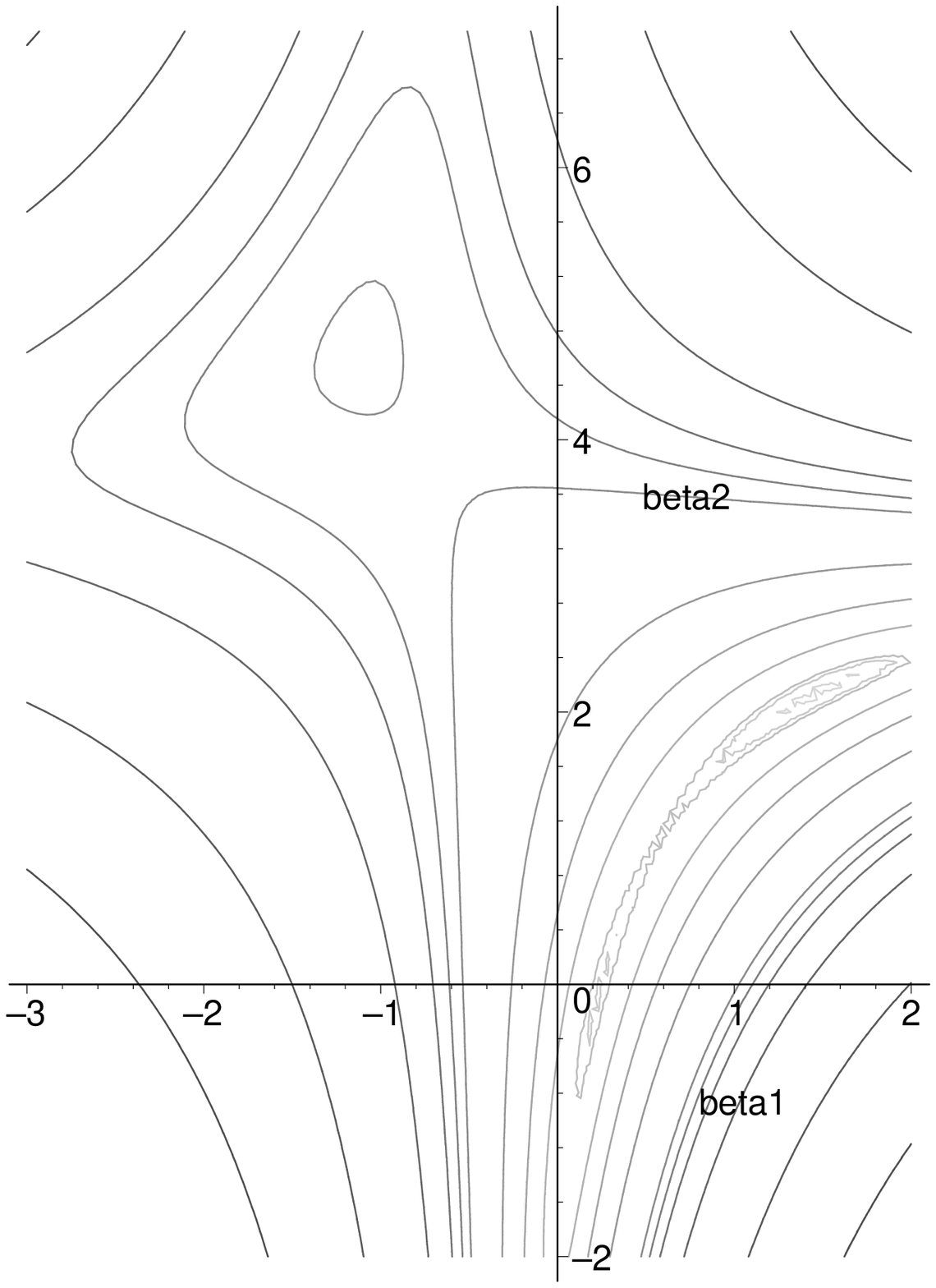} 
\end{center} 
\end{minipage}
\caption{\label{wireframe.3M} Three-dimensional plot and contour plot
  of profile log-likelihood function. 
%% Contour levels are at
%% -25, -23, {-21}, -20, -19.5, -19.05, -17, -15, {-13}, -11.2, -11,
%% and -10.7.    
}
\end{figure}

Table
\ref{tab:gensurdimdegree} states the dimension and degree of the
maximum likelihood ideal for seven examples including the
one from Section \ref{sec:revisit}.
\begin{table}[p]
  \begin{center}
    \begin{tabular}[h]{cccc} 
      $\CRfam$ & $\BCR$ & $\dim(I_G)$ & $\mathrm{degree}(I_G)$\\
      \hline
      \hline
      \raisebox{-0.25cm}{
      $\{(1,1),(2,2)\}$} & 
      \raisebox{-0.25cm}{
        $\left(\begin{smallmatrix}
        \beta_{11} & 0\\
        0 & \beta_{22}
      \end{smallmatrix}\right)$} & \raisebox{-0.25cm}{0} &
      \raisebox{-0.25cm}{5}\\[0.75cm] 
      $\{(1,1),(1,2),(2,3)\}$ & 
        $\left(\begin{smallmatrix}
        \beta_{11} & \beta_{12} & 0\\
        0 & 0 & \beta_{23}
      \end{smallmatrix}\right)$ & 0 & 9\\[0.5cm] 
    $\{(1,1),(2,2),(3,3)\}$ & 
    $\left(\begin{smallmatrix}
        \beta_{11} &  0 & 0\\
        0 & \beta_{22} & 0\\
        0 & 0 & \beta_{33}
      \end{smallmatrix}\right)$ & 0 & 29\\[0.5cm]
    $\{(1,1),(1,2),(1,3),(2,4)\}$ & 
    $\left(\begin{smallmatrix}
        \beta_{11} & \beta_{12} & \beta_{13} & 0\\
        0 & 0 & 0 & \beta_{24}
      \end{smallmatrix}\right)$ & 1 & 4\\[0.5cm]
    $\{(1,1),(1,2),(2,3),(2,4)\}$ & 
    $\left(\begin{smallmatrix}
        \beta_{11} & \beta_{12} & 0 & 0\\
        0 & 0 & \beta_{23} & \beta_{24}
      \end{smallmatrix}\right)$ & 1 & 8\\[0.5cm]
    $\{(1,1),(2,2),(3,3),(4,4)\}$ & 
    $\left(\begin{smallmatrix}
        \beta_{11} & 0 & 0 & 0\\
        0 & \beta_{22} & 0 & 0\\
        0 & 0 & \beta_{33} & 0\\
        0 & 0 & 0 & \beta_{44}
      \end{smallmatrix}\right)$ & 1 & 32\\[0.65cm]
    $\{(1,1),(2,2),(3,3),(4,4),(5,5)\}$ & 
    $\left(\begin{smallmatrix}
        \beta_{11} & 0 & 0 & 0 & 0\\
        0 & \beta_{22} & 0 & 0 & 0\\
        0 & 0 & \beta_{33} & 0 & 0\\
        0 & 0 & 0 & \beta_{44} & 0\\
        0 & 0 & 0& 0 & \beta_{55}
      \end{smallmatrix}\right)$ & 2 & 80\\[0.75cm]
        \hline  
    \end{tabular}
  \end{center}
  \smallskip
  \caption{Dimension and degree of maximum likelihood ideals.}
  \label{tab:gensurdimdegree}
\end{table}
For the models with zero-dimensional maximum likelihood ideal
$I_G$, we can find all stationary points of the likelihood function
by computations analogous 
to the ones demonstrated in Table \ref{tab:singularcode}.  
The likelihood functions of these models may be multimodal and
it would be interesting to find, for each model, reference data for which 
the cardinality of $V_\RRR(I_G)$ is large.  
For example, let $\CRfam=\{(1,1),(2,2)\}$ and choose
\begin{equation}
  \label{eq:xytrimod}
  \begin{split}
  X = &
  \left(\begin{array}{r@{\hspace{0.25cm}}r@{\hspace{0.25cm}}r@{\hspace{0.25cm}}r@{\hspace{0.25cm}}r} 
    -0.65 & -0.80 & \phantom{-}1.34  & -1.03 & -1.08 \\
    -0.04 & -1.18 & \phantom{-}1.98  & -2.42 & -3.75 
  \end{array}\right), \\
  Y = &
  \left(\begin{array}{r@{\hspace{0.25cm}}r@{\hspace{0.25cm}}r@{\hspace{0.25cm}}r@{\hspace{0.25cm}}r} 
    \phantom{-}0.14  & -0.73 & \phantom{-}1.40  & -2.29 & -3.30 \\
    \phantom{-}0.52 & -1.93 & 3.02 & -6.67 & -9.94
 \end{array}\right),
 \end{split}
\end{equation}
then the variety of the maximum likelihood ideal of $\NSUR$ is purely
real, i.e.\  
$V_\RRR(I_G)=V_\CCC(I_G)=5$. 
Figure \ref{wireframe.3M} shows a three-dimensional plot and a
contour plot of the profile log-likelihood function for these
observations. 
We conjecture 
that data with $V_\RRR(I_G)=V_\CCC(I_G)$ exist for all 
three models in Table \ref{tab:gensurdimdegree} that have
zero-dimensional maximum likelihood ideal. For 
models with  maximum likelihood ideal of dimension one or higher,
it is not clear whether $V_\RRR(I_G)=\infty$, i.e.\ a likelihood
function with an infinite number of stationary points, can occur with
non-zero probability.

\subsection{Submodels of seemingly unrelated regressions} 
\label{sec:submodels}

It is obvious that the algebraic approach developed in Section
\ref{sec:mlepoly} immediately carries over to
the submodels of SUR that are of interest
in testing equality of regression coefficients.  In the model $\NSUR$ with
$\CRfam=\{(1,1),(2,2)\}$, for example, we may be interested in testing
whether $\beta_{11}=\beta_{22}$.  If this is done using a likelihood ratio
test, then the likelihood function of the submodel in which
$\beta_{11}=\beta_{22}$ is imposed has to be maximized.  More precisely,
the submodel has the restricted parameter space
\begin{equation}
  \label{eq:subspace}
  \{ B\in \BCR\mid \beta_{11}=\beta_{22}\} \times \mathbb{P}.
\end{equation}
Table \ref{tab:submoddimdeg} lists similarly obtained submodels of the
models in Table \ref{tab:gensurdimdegree}, for which the maximum likelihood
ideal is zero-dimensional and the variety $V_\CCC(I_G)$ can be computed.
\begin{table}[t]
  \begin{center}
    \begin{tabular}[h]{cccc} 
      $\CRfam$ & Subspace of $\BCR$ & $\dim(I_G)$ & $\mathrm{degree}(I_G)$\\
      \hline
      \hline
      \raisebox{-0.25cm}{
      $\{(1,1),(2,2)\}$} & 
      \raisebox{-0.25cm}{
        $\left(\begin{smallmatrix}
        \beta_{11} & 0\\
        0 & \beta_{11}
      \end{smallmatrix}\right)$} & \raisebox{-0.25cm}{0} &
      \raisebox{-0.25cm}{3}\\[0.75cm] 
      $\{(1,1),(1,2),(2,3)\}$ & 
        $\left(\begin{smallmatrix}
        \beta_{11} & \beta_{12} & 0\\
        0 & 0 & \beta_{12}
      \end{smallmatrix}\right)$ & 0 & 7\\[0.5cm] 
    $\{(1,1),(2,2),(3,3)\}$ & 
    $\left(\begin{smallmatrix}
        \beta_{11} &  0 & 0\\
        0 & \beta_{11} & 0\\
        0 & 0 & \beta_{33}
      \end{smallmatrix}\right)$ & 0 & 11\\[0.5cm]
    $\{(1,1),(1,2),(1,3),(2,4)\}$ & 
    $\left(\begin{smallmatrix}
        \beta_{11} & \beta_{12} & \beta_{13} & 0\\
        0 & 0 & 0 & \beta_{13}
      \end{smallmatrix}\right)$ & 0 & 11\\[0.5cm]
    $\{(1,1),(1,2),(2,3),(2,4)\}$ & 
    $\left(\begin{smallmatrix}
        \beta_{11} & \beta_{12} & 0 & 0\\
        0 & 0 & \beta_{12} & \beta_{24}
      \end{smallmatrix}\right)$ & 0 & 23\\[0.5cm]
    $\{(1,1),(2,2),(3,3),(4,4)\}$ & 
    $\left(\begin{smallmatrix}
        \beta_{11} & 0 & 0 & 0\\
        0 & \beta_{11} & 0 & 0\\
        0 & 0 & \beta_{33} & 0\\
        0 & 0 & 0 & \beta_{44}
      \end{smallmatrix}\right)$ & 0 & 63\\[0.65cm]
    \hline
   \end{tabular}
  \end{center}
  \smallskip
  \caption{Dimension and degree of maximum likelihood ideals of submodels.}
  \label{tab:submoddimdeg}
\end{table}

It should also be noted that submodels of SUR need not inherit
unimodal likelihood 
functions from their parent model.  For example, the bivariate SUR
model $\NSUR$ with $\CRfam=\{(1,1),(2,1),(2,2)\}$ is monotone, 
i.e.\ the family $(C_r\mid r\in R)$ is totally ordered by inclusion, which
guarantees that the likelihood function has precisely one
stationary point corresponding to the global maximum
\citep{ap:ims,drton:2003}. However, the submodel induced by the
restriction $\beta_{11}=\beta_{21}$ can be reexpressed in the form of
the model studied in Section \ref{sec:revisit} by means of the linear
transformation that changes response $Y_2$ into $Y_2-Y_1$. Hence, the
submodel does not always have a unimodal likelihood function.

\section{Conclusion}
\label{sec:conclusion}

The presented algebraic approach to maximum likelihood estimation in
SUR permits us to 
compute all stationary points of the 
likelihood function if the maximum likelihood ideal is
zero-dimensional.   This is the case for three seemingly unrelated
regressions models considered in this paper (cf.\ Table
\ref{tab:gensurdimdegree}): (i) the previously studied model based on
$\CRfam^{(1)}=\{(1,1),(2,2)\}$, (ii) the model with
$\CRfam^{(2)}=\{(1,1),(1,2),(2,3)\}$, and (iii) the model with
$\CRfam^{(3)}=\{(1,1),(2,2),(3,3)\}$.   Additionally, 
interesting submodels of SUR may have a
zero-dimensional maximum likelihood ideal (cf.\ Table
\ref{tab:submoddimdeg}).  The computations in {\tt Singular} that
find all stationary points of the likelihood
functions of the models with zero-dimensional maximum likelihood ideal are
instantaneous for all but the model in Table \ref{tab:submoddimdeg}
that has a maximum likelihood ideal of degree 63.  
Thus we advocate the use of {\tt Singular} or similarly capable
software in statistical data analysis.  

In future work it would be interesting to find 
reference data sets leading to likelihood functions with a large
number of stationary points.  Moreover, the algebraic approach
presented herein could be 
combined with regression approaches \citep[e.g.][]{ap:ims,drton:2003} in
order to identify larger classes of SUR models for which all stationary
points of the likelihood function can be computed.  
Finally, it could be explored whether
methods for global minimization of polynomials
\citep{parrilo:2001} can be used to find the global maximum of SUR
likelihood functions.

\begin{ack}
I would like to thank Michael Perlman, Thomas Richardson, and Bernd
Sturmfels for their help and support, and an anonymous referee for helpful
comments on the presentation.
\end{ack}

\end{document}